\documentclass[12pt,twoside,reqno]{amsart}
\usepackage{amssymb, latexsym}
\usepackage[dvipsnames]{xcolor}

\newcommand{\beqa}{\begin{eqnarray*}}
\newcommand{\eeqa}{\end{eqnarray*}}
\newcommand{\beqn}{\begin{eqnarray}}
\newcommand{\eeqn}{\end{eqnarray}}

\newcommand{\mcH}{\mathcal H}

\newcommand{\mcL}{\mathcal L}

\newcounter{cnt1}
\newcounter{cnt2}
\newcounter{cnt3}
\newcommand{\blr}{\begin{list}{$($\roman{cnt1}$)$}
 {\usecounter{cnt1} \setlength{\topsep}{0pt}
 \setlength{\itemsep}{0pt}}}
\newcommand{\bla}{\begin{list}{$($\alph{cnt2}$)$}
 {\usecounter{cnt2} \setlength{\topsep}{0pt}
 \setlength{\itemsep}{0pt}}}
\newcommand{\bln}{\begin{list}{$($\arabic{cnt3}$)$}
 {\usecounter{cnt3} \setlength{\topsep}{0pt}
 \setlength{\itemsep}{0pt}}}
\newcommand{\el}{\end{list}}
\newtheorem{thm}{Theorem}[section]

\newtheorem{cor}[thm]{Corollary}

\newtheorem{Def}[thm]{Definition}

\newtheorem{rem}[thm]{Remark}
\newcommand{\Rem}{\begin{rem} \rm}
\newcommand{\bdfn}{\begin{Def} \rm}
\newcommand{\edfn}{\end{Def}}

\newcommand{\ba}{\begin{array}}
\newcommand{\ea}{\end{array}}

\begin{document}
\begin{center}\large{{\bf{ Modular convergence in $H$-Orlicz spaces of Banach valued functions}}} 
 
  {\tiny{ Hemanta Kalita$^{1, \ast}$  
  and  Bipan Hazarika$^{2}$\\

$^1$Department of Mathematics, Assam Don Bosco University, Sonapur, Guwahati 782402, Assam, India\\
\noindent$^{2}$ Department of Mathematics, Gauhati University, Guwahati  781014, Assam, India\\Email:hemanta.kalita@dbuniversity.ac.in;hemanta30kalita@gmail.com; 
bh\_rgu@yahoo.co.in; bh\_gu@gauhati.ac.in}}

\end{center}
\title{}
\author{}
\thanks{$^{\ast}$The corresponding author}

\begin{abstract} In this article we develop the theory of $H$-Orlicz
space generated by generalised Young function. Modular convergence of $H$-Orlicz space for the case of vector-valued functions  and norm convergence in $\mcH^\theta(X, \overline{\mu})$ where $X$ is any Banach space are discussed. Relationships of modular convergence and norm convergence of  $H$-Orlicz spaces are discussed.
\\
\noindent{\footnotesize {\bf{Keywords and phrases:}}} H-Orlicz space, Modular convergence, Norm convergence\\
{\footnotesize {\bf{AMS subject classification \textrm{(2020)}:}}} 26A39, 46B03, 46B20, 46B25.
\end{abstract}
\maketitle


\pagestyle{myheadings}
\markboth{\rightline {\scriptsize   HK}}
        {\leftline{\scriptsize  }}

\maketitle

\section{Introduction and preliminaries}
Z.W. Birnbaum and W. Orlicz to proposed a generalized space of $L^p,$ later on it was known as Orlicz space. This space was later developed by Orlicz himself. The fundamental properties
of Orlicz space with Lebesgue measure found in \cite{MA} also (see \cite{Mali}). H. Nakano \cite{Nakano} introduced the concept of modular Orlicz space also (see \cite{MO}) for modular space.  M.S. Skaff \cite{MS} developed generalised N-function. Gereralised N-function are generalization of the variable N-function ( see \cite{Por}) and non decreasing N-function (see \cite{Wang}). M.S. Skaff \cite{MS1} discussed vector valued Orlicz spaces with generalised N-function.  A. Kozek \cite{Koz} studied  Orlicz spaces of functions with values in Banach spaces. In his work, Kozek discussed Orlicz spaces from the modular spaces from point of view of generalised N-function. A. Kaminska and H. Hudzik \cite{KH} discussed the necessary condition of equality of the  modular convergence and norm convergence in Orlicz spaces. In recent time, modular convergence theorem is a well known concept for many researcher in different areas of mathematics. Carlo Bardaro and Gianluca Vint \cite{Carlo} studied modular convergence theorem for certain nonlinear integral
operators with homogeneous kernel,  modular convergence in fractional Musielak-Orlicz Spaces (see \cite{Carlo1}). Youssef Ahmida et.al.,  discuss density of smooth functions on Musielak-Orlicz-Sobolev spaces (see \cite{You}). Hazarika and Kalita \cite{BH} introdced H-Orlicz spaces with non absolute integrable functions in particluar Henstock-Kurzweil integrable functions. The significance of H-Orlicz space is $C_0^\infty$ is dense in H-Orlicz space which is not generally true in the case of Orlicz spaces. On the other side about Henstock-Kurzweil integral in 1957, Jaroslav Kurzweil discussed about a new integral in one of his publication, while unaware of the work of Kurzweil, Ralph Henstock published an article on integration theory in which he discussed the same integration as J. Kurzweil. This new integral can integrate a substantial type of functions compare to  the Riemann or Lebesgue integral. In the honors of these mathematicians, nowadays this integral is called Henstock-Kurzweil integral in brief $HK$-integral. Measure theory is not essential in the definition of HK-integral. In quantum theory and nonlinear analysis, HK-integrals are aid for highly oscillatory functions to integrate.  Moreover, HK integrability encloses 
improper integrals (see \cite{AB1,ABH,SAR2,BH}). Major drawback of Henstock-Kurzweil integrable function space is not complete with Alexiewicz norm. In this article in Section 2, we introduce a new norm which is equivalent to Alexiewicz norm but Henstock-Kurzweil integrable function space become complete with this particular norm. In Section 3, we extend the theory of H-Orlicz spaces with vector functions from the point of view of generalised Young-function. Finally in Section 3, we discuss various relationship of modular convergence as well as norm convergence of Orlicz spaces and H-Orlicz spaces.
\\Let $(\mathcal{J}, \mathfrak{Q}, \overline{\mu})$ be a measure space, where $\mathcal{J}$ is an abstract set, $\mathfrak{Q}$ is a $\sigma-$algebra of subsets of the set $\mathcal{J},~\overline{\mu}$ is a $\sigma-$finite, positive, complete measure on $\mathfrak{Q}.~X$ is a Banach space. Let $M_X$ be a set of all $\overline{\mu}-$measurable functions $\xi:~\xi^i(t)~( t \in \mathcal{J},~i=1,2,..)$ where $f^i(t)$ are real valued function defined on $\mathcal{J}$ with values in $X.$
\begin{Def}\label{def,}
\cite[Def 2.1]{KH}
A function $\theta: X \times \mathcal{J} \to [0, \infty]$ is called an $N-$function if it satisfies the following conditions:\\ There exists a set $\mathcal{J}_1 \in \mathfrak{Q},~\overline{\mu}( \mathfrak{Q} \setminus \mathcal{J}_1) =0,$ such that 
\begin{enumerate}
\item[a)] $\theta(.,.)$ is $B \times \mathfrak{Q}$ measurable, where $B$ denotes the $\sigma-$algebra of Borel subsets of $X;$
\item[b)] $\theta(., t)$ is lower semi continuous on $X$ for every $t \in \mathcal{J}_1;$
\item[c)] $\theta(., t)$ is convex for every $t \in \mathcal{J}_1$
\item[d)] $\theta(0,t)=0$ and $\theta(\xi(t),t)= \theta(-\xi(t),t)$ for every $\xi \in M_X,~t \in \mathcal{J}_1;$
\item [e)] there exist $\overline{\mu}-$measurable functions $\alpha: \mathcal{J}_1 \to (0, \infty)$ and $\lambda: \mathcal{J}_1 \to (0, \infty)$ such that $||\xi(t)|| \geq \lambda(t)$ implies $\theta(\xi(t),t) \geq \alpha(t)$ for every $t \in \mathcal{J}_1$
\item[f)] there exist $\overline{\mu}-$measurable functions $\rho: \mathcal{J}_1 \to (0, \infty)$ and $\rho_0: \mathcal{J}_1 \to (0, \infty)$ such that $||\xi(t)|| \leq \rho(t)$ implies $\theta(\xi(t),t) \leq \rho_0(t)$ for every $t \in \mathcal{J}_1.$
\end{enumerate}
\end{Def}
If $\theta$ fulfills conditions $(a)-(d),$ then $\theta$ is called $N^{''}-$function; if it satisfies condition $(a)-(e),$ it is called $N^{'}-$function. Recall that an $N^{'}-$ function $\theta$ is called an $N-$function if $X$ is separable, reflexive Banach space if for almost $t \in \mathcal{J},$
\begin{eqnarray*}
\theta(x,t) \leq \xi(t) < \infty~for~||x|| \leq \rho(t),~\xi(t),~\rho(t) \in (0, \infty)~holds.
\end{eqnarray*}
 A function $\xi: \mathcal{J} \to X$ is $M$-measurable if there is a sequence of simple functions from $S_X(\mathfrak{Q})$ converges to $\xi~(M-$a.e), where $S_X(\mathfrak{Q})$ is a set of all $X-$valued simple functions.\\ $M$ denotes a $\sigma-$bounded family of positive measures defined on $\mathfrak{Q}.$ This means for each $E \in \mathfrak{Q},$ there exists a pairwise disjoint collections $\{E_i\}_{i=1}^{\infty},~E_i \in \mathfrak{Q}$ such that $E =  \bigcup_{i=1}^{\infty} E_i$ and $S_X(E)$ is a set of all $X$ valued simple functions.
\begin{Def}
\cite{MO} Given a linear space $X,$ a functional $\rho(.)$ defined on $X$ with values $-\infty< \rho(x) \leq \infty$ is called a modular if the following conditions hold:
\begin{enumerate}
\item[(a)] $\rho(x)=0$ if and only if $x=0$
\item[(b)] $\rho(-x)= \rho(x),$
\item[(c)] $\rho( \alpha x + \beta y) \leq \rho(x)+\rho(y)$ for every $\alpha, \beta\geq0,~\alpha+\beta=1.$
\end{enumerate}
\end{Def}
\begin{Def}

\cite{MO}
A sequence $(x_n) \subset X$ is said to be modular convergent to $x \in X$ if there exists a number $\alpha>0 $(depending on the sequence $(x_n)$) such that $\rho(\alpha(x_n -x)) \to 0$ as $n \to \infty.$
\end{Def}
\begin{Def}
\cite[Definition 2.2]{MS1}
Let $\theta(t, \xi(t))$ be a generalised Young function and let $L^\theta(X, \overline{\mu})$ be its associated Orlicz class. The closure $\mcL^\theta(X, \overline{\mu})$ of $L^\theta(X, \overline{\mu})$ under positive scalar multiplication is a vector valued Orlicz space.  
\end{Def}
The space $\mcL^\theta(X, \overline{\mu})$ is a Banach space with Luxemberg type norm 
\begin{eqnarray*}
||\xi||_{\mcL^\theta}= \inf \bigg\{\alpha>0:~(L)\int_{\mathcal{J}} \theta\bigg(t, \frac{\xi(t)}{\alpha}\bigg)d \overline{\mu} \leq 1 \bigg\}.
\end{eqnarray*}
Recalling convergence of a sequence of function $\xi_n$ in the normed vector space $\mcL^\theta(X, \overline{\mu})$ to $\xi$ in following way: 
\begin{eqnarray*}
\xi_n \to \xi~if~||\xi_n-\xi||_{\mcL^\theta} \to 0~as~n \to \infty.
\end{eqnarray*}
\section{Henstock-Kurzweil integrable function spaces $H(X, \overline{\mu})$}
In this section of the paper we discuss about the Henstock-Kurzweil integrable function space. We have leave off the basic results of the space. We discuss only those results that are closely connected for next section.\\
A finite collection $D=\{D_i:~i=1,2,..,n\}$ of mutually disjoint elements of $\mathfrak{Q}$ is called sub-partition of $\mathcal{J}$ if $\mathcal{J}= \bigcup_{i=1}^{n}D_i.$ In this case $D=\{D_i:~i=1,2,..,n\}$ is called partitions of $\mathcal{J}.$\\ Let $\widehat{D}=\{(D_i, d_i):~i=1,2,..,n\}$ and $d_i \in D_i.$ Clearly the norm of the sub-partition $D=\{D_i\}_{i=1}^{n}$ is $||D||= \sup\{\overline{\mu}(D_i):~i=1,2,..,n\}.$ We assume the collection of all tagged sub-partitions of $\mathcal{J}$ will be $\mathfrak{Y}.$ If $\widehat{D_1},~\widehat{D_2} \in \mathfrak{Y}$ then the sub partition
  \begin{eqnarray*}
  D_1 \vee D_2= (D \setminus \cup D_2) \cup \{A \cap B:~A \in D_1,~B \in D_2\}.
  \end{eqnarray*}
  Clearly $(\mathfrak{Y}, \gg)$ is a directed set. Let $\xi: \mathcal{J} \to X$ and $\widehat{D}= \{(\widehat{D_1}, d_i):~i=1,2,..,n\} \in \mathfrak{Y}$ then we construct the Riemann sum of $\xi$ corresponding to $\widehat{D}$ as $$S(\xi, \widehat{D})= \sum_{i=1}^{n}\xi(d_i) \overline{\mu}(D_i).$$  As $(\mathfrak{Y}, \gg)$ is directed set, the function $\widehat{D} =\{(D_i,d_i):~i=1,2,..,n\} \to S(\xi, \widehat{D})$ defines a net $S: \mathfrak{Y} \to {X}.$
  \begin{Def}
  A function $\xi: \mathcal{P} \subset \mathfrak{Q} \to X$ is said to be Henstock-Kurzweil integrable on the set $\mathcal{P} \in \mathfrak{Q}$ if there is an element $I_{\mathcal{P}}= \int_{\mathcal{P}}f d\overline{\mu}$ such that for every $\epsilon>0$ there exists $\widehat{D_0} \in \mathfrak{Y}$ and $\widehat{D} \gg \widehat{D_0}$ in $\mathfrak{Y},$ then $$||S(\xi, \widehat{D}) - \int_{\mathfrak{Q}}\xi d \overline{\mu}|| \leq \epsilon.$$
  \end{Def}
 We call net limits of  the integral $\xi: \mathfrak{Q} \to X$ as $$\int_{\mathfrak{Q}}\xi d\overline{\mu} = \lim\big(S(\xi, \widehat{D})_{D \in (\mathfrak{Y}, \gg )}\big).$$ This integrals are unique, linear also for  all $X-$valued Henstock integrable. Let $HK(\mathfrak{Q}, \overline{\mu},X) $ or in brief $HK(X, \overline{\mu})$ the set of all ${X}$-valued $\overline{\mu}-$ Henstock-Kurzweil integrable functions. 
 \begin{Def}
 For each $\xi:\mathfrak{Q} \to X$ 
 \begin{eqnarray}\label{bh}
 ||\xi||_{\mathfrak{Y}}= \sup \{||S(\xi, \widehat{D})||:~\widehat{D} \in \mathfrak{Y}\}
 \end{eqnarray} 
 \begin{thm}
  The expression (\ref{bh}) is a norm on $HK(X, \overline{\mu}).$ 
 \end{thm}
 \begin{proof}
 For any vector valued function $\xi_1, \xi_2 \in HK(X, \overline{\mu}),$ we have the non negativity as follows:
  \begin{align*}
 {\textit{(i)}}~||\xi_1||_{\mathfrak{Y}} &= \sup \bigg\{||S(\xi_1, \widehat{D})||~:~\widehat{D} \in \mathfrak{Y} \bigg\} \\& =\sup \bigg\{||\sum_{i=1}^{n}\xi_1(d_i)\overline{\mu}(D_i)||~:~\widehat{D}=(D_i, d_i),~i=1,2,..\bigg\}\\& \geq 0
 \end{align*}
 Therefore, $||\xi_1||_{\mathfrak{Y}} \geq 0.$\\
 (ii) Now, $||\xi_1||_{\mathfrak{Y}}=0~$ if and only if 
 \begin{align*}
 ||\xi_1||_{\mathfrak{Y}} &= \sup \bigg\{||S(\xi_1, \widehat{D})||~:~\widehat{D} \in \mathfrak{Y} \bigg\} \\& =\sup \bigg\{||\sum_{i=1}^{n}\xi_1(d_i)\overline{\mu}(D_i)||~:~\widehat{D}=(D_i, d_i),~i=1,2,..\bigg\}\\& = \sum_{i=1}^{n}\xi_1(d_i)\overline{\mu}(D_i) \\& =0
 \end{align*}
 Therefore, $||\xi_1||_{\mathfrak{Y}}=0~$ if and only if  $\xi_1 =0.$\\
 (iii) ~for a scalar $\alpha,$ we have
 \begin{align*}
 ||\alpha \xi_1||_{\mathfrak{Y}} &= \sup \bigg\{||S(\alpha \xi_1, \widehat{D})||~:~\widehat{D} \in \mathfrak{Y} \bigg\} \\&= \sup \bigg\{||\alpha S(\xi_1, \widehat{D})||~:~\widehat{D} \in \mathfrak{Y} \bigg\}~{\textit{\cite[page 199]{Bartle}}} \\&= |\alpha| \sup \bigg\{||S(\xi_1, \widehat{D})||~:~\widehat{D} \in \mathfrak{Y} \bigg\}~\\& =|\alpha|||\xi_1||_{\mathfrak{Y}}
 \end{align*}
 (iv) for triangle inequality,
 \begin{align*}
 ||\xi_1+\xi_2||_{\mathfrak{Y}} &= \sup \bigg\{||S(\xi_1+\xi_2, \widehat{D})||~:~\widehat{D} \in \mathfrak{Y} \bigg\} \\& = \sup \bigg\{||S(\xi_1, \widehat{D})+ S(\xi_2, \widehat{D})||~:~\widehat{D} \in \mathfrak{Y} \bigg\}~{\textit{\cite[page 199]{Bartle}}}\\& \leq \sup \bigg\{||S(\xi_1, \widehat{D})||~:~\widehat{D} \in \mathfrak{Y} \bigg\}+\sup \bigg\{||S(\xi_2, \widehat{D})||~:~\widehat{D} \in \mathfrak{Y} \bigg\} \\& \leq ||\xi_1||_{\mathfrak{Y}} +||\xi_2||_{\mathfrak{Y}}
 \end{align*}
 \end{proof}
 Then clearly $||\xi||_{\mathfrak{Y}} < \infty,$ and the norm (\ref{bh}) is equivalent to the Alexiecz norm for the Henstock integrable function spaces.
 \end{Def}
 \begin{thm}\label{4rd}
 Let $\xi: \mathfrak{Q} \to X.$ If $\xi=0~\overline{\mu}-$a.e. on $\mathfrak{Q},$ then $\xi \in HK(X, \overline{\mu})$ and $\int_{\mathfrak{Q}}\xi d \overline{\mu}=0.$
 \end{thm}
 \begin{proof}
 Let $\mathcal{B}=\{z \in \mathfrak{Q}:~\xi(z) \neq 0\}.$ For each positive $n,$ assume 
 \begin{eqnarray*}
 \mathcal{B}_n = \{z \in \mathcal{B}:~n-1 \leq ||\xi(z)||<n\}
 \end{eqnarray*}
 Choosing an open set $\mathcal{O}_n$ for each $n$ such that $\mathcal{B}_n \subseteq \mathcal{O}_n$ and $\overline{\mu}(\mathcal{O}_n) < \frac{\epsilon}{n. 2^n}.$ For a positive function \[
\xi(x)  = \left\{ {\begin{array}{*{20}c}
   1~if~x \in \mathfrak{Q} -\mathcal{B}  \\
   e~if ~x \in \mathcal{B}_n  \\

 \end{array} } \right.
\] where $e(z, \mathcal{B}_n)= \inf\{|y-z|:~y \in \mathcal{B}_n\}$ then $e(z, \mathcal{B}_n)>0$ if $x \notin \mathcal{B}_n$ and $\mathcal{B}_n$ is closed. If $D_n$ be the subset of tagged partition $D$ of $\mathfrak{Q}$ that is sub-ordinate to $\xi$ in $\mathcal{B}_n.$\\Now,
\begin{align*}
||S(\xi, \widehat{D})|| &\leq \sum_{n=1}^{\infty}||S(\xi, \widehat{D})||\\&< \sum_{n=1}^{\infty}n \nu(\mathcal{O}_n)\\&< \sum_{n=1}^{\infty} \epsilon. 2^{-n}\\&= \epsilon
\end{align*}
So, $\xi \in HK(X, \overline{\mu})$ and $\int_{\mathfrak{Q}}\xi=0.$
 \end{proof}
 \begin{cor}
 Let $\xi,g \in HK(X, \overline{\mu})$ be $\overline{\mu}-$essentially equal then $\int_{\mathfrak{Q}}\xi d\overline{\mu}= \int_{\mathfrak{Q}}g d\overline{\mu}.$
 \end{cor}
 \begin{thm}
 Let $(\mathfrak{Q}, \mathcal{J}, \overline{\mu}) $ be a finite measure, $X$ is a Banach space then $\big(HK( X, \overline{\mu}), ||.||_{\mathfrak{Y}}\big)$ is complete.
 \end{thm}
 \begin{proof}
 Let us assume $(\xi_n)$ be a Cauchy sequence in $HK(X, \overline{\mu}),$ then for each $\epsilon>0$ we can find a natural number $N$ such that $m, n >N$ such that 
 \begin{eqnarray*}
 \sup\limits_{\widehat{D} \in \mathfrak{Y}}||S(\xi_m-\xi_n, \widehat{D})|| \leq \epsilon.
 \end{eqnarray*}
 If $\overline{w} \in \mathfrak{Q},$ under the assumption of sub-partition $\widehat{D}=\{(\mathfrak{Q}, \overline{w})\},$ we get $||\xi_n(\overline{w})- \xi_m(\overline{w})|| \leq \overline{\mu}\big( \mathfrak{Q}\big) \epsilon.$ From our hypothesis $\epsilon>0,~\overline{\mu}(\mathfrak{Q}) < \infty,~$ we find $\big(\xi_n(\overline{w})\big)$ is Cauchy sequence of $X.$ Again, since $X$ is a Banach space, it is easy to define $\overline{w} \to \xi(\overline{w})= \lim\limits_{n \to \infty}\xi_n(\overline{w}).$\\ Using the concept of the Riemann sum over sub-partitions $\widehat{D}_1$ and $\widehat{D}_2$ in $\mathfrak{Y}$ with $\widehat{D}= \widehat{D}_1 \cup \widehat{D}_2$ we get the following:
 \begin{align*}
& ||\int_{\mathfrak{Q}}\xi_n d \overline{\mu}- \int_{\mathfrak{Q}}\xi_m d \overline{\mu}||\\  &\leq ||\int_{\mathfrak{Q}}\xi_n d \overline{\mu}- S(\xi_n, \widehat{D})||  +||S(\xi_n, \widehat{D})-S(\xi_m, \widehat{D})||+||\int_{\mathfrak{Q}}\xi_m d \overline{\mu}- S(\xi_n, \widehat{D})||\\& < 3 \epsilon.
 \end{align*}
 This gives, for $N \in \mathbb{N}, ~||\int_{\mathfrak{Q}}\xi_n d \overline{\mu}-\int_{\mathfrak{Q}}\xi_m d \overline{\mu}|| < \epsilon $ for $m,n~\geq N.$ So, we can conclude $\int_{\mathfrak{Q}}\xi_n d \overline{\mu}$ is Cauchy sequence in $X.$ Say $\int_{\mathfrak{Q}}\xi_n d \overline{\mu}$ converges to $x \in X.$\\ Again,
 \begin{align*}
& ||S(\xi, \widehat{D})-x||\\ & \leq ||S(\xi, \widehat{D})- S(\xi_n, \widehat{D})||+||S(\xi_n, \widehat{D})-S(\xi_m, \widehat{D})||\\ &+||S(\xi_m, \widehat{D})-\int_{\mathfrak{Q}}\xi_m d \overline{\mu}||+||\int_{\mathfrak{Q}}\xi_m d \overline{\mu}-x|| \\& < 4 \epsilon.
 \end{align*}
 So, $\xi \in HK(X, \overline{\mu})$ and $x= \int_{\mathfrak{Q}} \xi d \overline{\mu}.$
 \end{proof}
 \begin{thm}
 If $X$  is weakly sequentially complete, then $HK\big(X, \overline{\mu}\big)$  is also sequentially complete.
 \end{thm}
 \begin{proof}
 The proof is similar as \cite[Theorem 2.2]{Fernando}. 
 \end{proof}
 \begin{thm}
 Let $\xi: \mathbb{J} \subseteq \mathcal{J} \to X$ be in $HK(X,\overline{\mu}),$ \\then $\theta\big(t,\xi(t) \big) \in HK(X, \overline{\mu})$ for all $t \in \mathbb{J}.$
 \end{thm}
 \begin{proof}
  Let $ \xi:\mathbb{J} \subseteq \mathcal{J}  \to X $ be   measurable (integrable) function, then for all $\epsilon > 0 $ there exists a $\delta: \mathbb{J} \to (0, +\infty) $ such that  $$\left|\left|\sum\limits_{(\mathbb{J},t) \in \pi}\sum\limits_{(\mathbb{J}^{'}, t^{'}) \in \pi^{'}} [\xi(t)-\xi(t^{'})]\overline{\mu}(\mathbb{J} \cap \mathbb{I}^{'})\right|\right|_X < \epsilon $$ for all partitions $\pi$ and $\pi^{'} $ of $\mathbb{J}$ finer than $\delta.$\\
  As $\theta $ is Young function. So, $\theta(t, \xi(t)) \to \infty $ as $ t\to \infty.$ Our claim is $\theta(t, \xi(t)) $ is  Henstock-Kurzweil integrable.\\
   Since, Young function by definition, is an extended real Borel function. So, $\theta(t, \xi(t)) $ is measurable. If $\pi$ and $\pi^{'} $ are both partitions of the same interval of  $\mathbb{J},$ then for any subinterval $ \mathbb{J}_0$ of $ \mathbb{J}_1$  we can write 
   \begin{align*} 
   &\overline{\mu}(\mathbb{J}) = \sum\limits_{(\mathbb{J}^{'},t^{'}) \in \pi^{'}}\overline{\mu}(\mathbb{J} \cap \mathbb{J}^{'}).\end{align*} 
   So, \begin{align*}
   &~~~~\sum\limits_{(\mathbb{J},t) \in \pi} \theta(t, \xi(t))(t) \overline{\mu}(\mathbb{J}) =\sum\limits_{(\mathbb{J},t) \in \pi}\sum\limits_{(\mathbb{J}^{'},t^{'}) \in \pi^{'}} \theta(t,\xi(t))\overline{ \mu}(\mathbb{J} \cap \mathbb{J}^{'})\\
&\mbox{i.e.,}~\left|\left|\sum\limits_{(\mathbb{J},t) \in \pi} \theta(t, \xi(t))\overline{ \mu}(\mathbb{J})  - \sum\limits_{(\mathbb{J}^{'},t^{'}) \in \pi^{'}} \theta(t,\xi(t))\overline{ \mu}(\mathbb{J} )\right|\right|_X < \epsilon.\end{align*}
   Thus, $ \theta(t,\xi(t))$ is in $HK(X, \overline{\mu}).$
  \end{proof}
  \begin{cor}
For all functions $\xi \in M_X$ for which there exists a constant $k>0$ such that 
\begin{eqnarray*}
\rho(k\xi)=  \int_{\mathcal{J}}\theta(t, k \xi(t)) d \overline{\mu} \in HK(X, \overline{\mu}).
\end{eqnarray*}
\end{cor}
\section{H-Orlicz spaces $\mcH^\theta(X, \overline{\mu})$}
 In this section, we initiate to study of H-Orlicz spaces associate with modular of functions. Let $M_X$ be the set of all real-valued (or complex-valued), $\mathfrak{Q}-$measurable and finite $\overline{\mu}$-a.e. functions on $\mathcal{J},$ with equality $\overline{\mu}$-a.e. Clearly $\theta(t, \xi(t))$ is $\mathfrak{Q}-$measurable function of $t \in \mathcal{J}$ for every $f \in M_X$ we define
\begin{eqnarray}\label{feb18}
\rho(\xi)= (H) \int_{\mathcal{J}}\theta(t, \xi(t))d \overline{\mu}.
\end{eqnarray}
Clearly,
$\rho(x)$ of the equation (\ref{feb18}) is a modular in $M_X.$

 We define the $H-$Orlicz class  as follows:
\begin{eqnarray*}
H^\theta(X, \overline{\mu})=\{ \xi: \mathcal{J} \to X~measurable:~\int_{\mathcal{J}}\theta\big(t, k \xi(t)\big)d\overline{\mu} \in H(X, \overline{\mu}),~for~some~k>0\}
\end{eqnarray*}
It is very straight forward that:
\begin{eqnarray*}
H^\theta(X, \overline{\mu})=\{ \xi: \mathcal{J} \to X~measurable:~\int_{\mathcal{J}}\theta\big(t, k \xi(t)\big)d\overline{\mu} \in H(X, \overline{\mu})~\} \to 0 ~as~ k \to 0+.
\end{eqnarray*}

\begin{thm}
The $H$-Orlicz class $H^\theta(X, \overline{\mu})$ is a convex set of functions. That is for given $\alpha, \beta \geq 0,~\alpha +\beta=1,~\rho(\alpha \xi + \beta g) \leq \alpha \rho(\xi)+\beta \rho(g).$
\end{thm}
\begin{proof}
The proof is similar as \cite[Theorem 2.1]{BH}.
\end{proof}
\begin{thm}
The $H$-Orlicz class $H^\theta(X, \overline{\mu})$ is linear if and only if $H^\theta(X, \overline{\mu})$ is closed under positive scalar multiplication.
\end{thm} 
\begin{proof}
Let $H^\theta(X, \overline{\mu})$ be a  linear space then, clearly it is closed under positive scalar multiplication.\\ Conversely, assume if $H^\theta(X, \overline{\mu})$ is closed under positive scalar multiplication. We will prove $H^\theta(X, \overline{\mu})$ is linear space. From the definition of $\theta(t, \xi(t)),$ we find $\theta(t, \xi(t))=\theta(t, -\xi(t)),$ this means $-\xi \in H^\theta(X, \overline{\mu}).$ Assume $ \alpha,~ \beta>0$ are real numbers and $\xi, g \in H^\theta(X, \overline{\mu}).$\\{\bf Case I:} Let $\alpha,~ \beta >0$ then for each $t \in \mathcal{J}$ and  convexity of $\theta\big(t, \xi(t)\big)$ we have
\begin{eqnarray}\label{2nd}
\theta\bigg(t, \frac{\alpha \xi(t)+ \beta g(t)}{\alpha +\beta}\bigg)= \theta\bigg(t, \frac{|\alpha| \xi(t)}{|\alpha|+|\beta|}\bigg) +\theta\bigg(t, \frac{|\beta| g(t)}{|\alpha|+|\beta|}\bigg)
\end{eqnarray}
 By the assumption $\alpha +\beta >0.$ The right sides of the equation (\ref{2nd}) is in $H(X, \overline{\mu})$ so we can conclude $\alpha \xi +\beta g \in H^\theta(X, \overline{\mu}).$\\{\bf Case II:} If $\alpha \beta <0$ with $\alpha <0< \beta$ then
\begin{eqnarray}\label{3nd}
\theta\bigg(t, \frac{|\alpha| (- \xi(t))+ \beta g(t)}{|\alpha| +\beta}\bigg)= \theta\bigg(t, \frac{|\alpha|(- \xi(t))}{|\alpha|+\beta}\bigg) +\theta\bigg(t, \frac{\beta g(t)}{|\alpha|+\beta}\bigg)
\end{eqnarray}
By the assumption $|\alpha| +\beta >0.$ The right sides of the equation (\ref{3nd}) is in $H(X, \overline{\mu})$ so we have $\alpha \xi +\beta g \in H^\theta(X, \overline{\mu}).$
\end{proof}
We define vector valued $H$-Orlicz space as below
\begin{Def}
The closure of  $H$-Orlicz class $H^\theta(X, \overline{\mu})$  under positive scalar multiplication of a vector valued  generalised Young function $\theta(t, \xi(t))$ will be called vector valued $H$-Orlicz space, denoted as $\mcH^\theta(X, \overline{\mu}).$ That is:
\begin{eqnarray*}
\mcH^\theta(X, \overline{\mu})= \big\{\xi \in H^\theta(X, \overline{\mu}):~ c \xi \in H^\theta(X, \overline{\mu}) \big\}
\end{eqnarray*}
\end{Def}

We define the norm of $H^\theta(X, \overline{\mu})$ as follows:
\begin{eqnarray}\label{3rd-}
||\xi||_{(X, \overline{\mu})}= \inf \bigg\{k>0:~(H)\int_{\mathcal{J}} \theta\bigg(t, \frac{\xi(t)}{k}\bigg) d \overline{\mu} \leq 1 \bigg\}
\end{eqnarray}
It is very clear that $\bigg(\mcH^\theta\big(X, \overline{\mu}\big), ||.||_{(X, \overline{\mu})}\bigg)$ is a Banach spaces with the norm (\ref{3rd-}).

\begin{thm}
 The classical Orlicz space $\mcL^\theta\big(X, \overline{\mu}\big)$ is a dense subspace of $\mcH^\theta\big(X, \overline{\mu}\big)$  as continuous dense embeddings. That is, $\mathcal{L}^\theta(X, \overline{\mu}) \hookrightarrow  \mcH^\theta\big(X, \overline{\mu}\big)$ is continuous dense embeddings.
\end{thm}
 \begin{proof}
    Let $h \in \mathcal{L}^\theta(X, \overline{\mu}).$ Then $h \in L^1(X, \overline{\mu})$ with $||f||_{L^\theta}< \infty.$ Then for some $k>0,$   we have
    \begin{align*}
  \inf\left\{ (H)\int_{\mathcal{J}} \theta\left(t,\dfrac{h(t)}{k}\right)d\overline{\mu} \right\} &\leq \inf\left\{(L)\int_{\mathcal{J}}\theta\left(t, \dfrac{h(t)}{k}\right)d\overline{\mu} \right\}\\& \leq 1.
   \end{align*}
   So, for some $k>0,~\inf\left\{k>0:~(L)\int_{\mathcal{J}} \theta\left(t, \frac{h(t)}{k}\right)d\overline{\mu} \leq 1\right\},$ we get the following $$\inf\left\{k>0:~(H)\int_{\mathcal{J}} \theta\left(t, \dfrac{h(t)}{k}\right)d\overline{\mu} \leq 1\right\}.$$  Hence $h \in \mathcal{H}^\theta(X, \overline{\mu})$ with $||h||_{(X, \overline{\mu})} \leq ||h||_{\mathcal{L}}.$ Hence the proof.
    \end{proof}
\begin{thm}
	Suppose $\overline{\mu}(X)<\infty$ and $\overline{\mu}$ is bounded, then $\mathcal{H}^\theta(X, \overline{\mu}) \hookrightarrow  L^1(X, \overline{\mu})$ is continuous.
\end{thm}
\begin{cor}
Suppose $\overline{\mu}(X)<\infty$ and $\overline{\mu}$ is bounded, then $\mathcal{H}^\theta(X, \overline{\mu}) \hookrightarrow  HK(X, \overline{\mu})$ is continuous.
\end{cor}
\section{ Modular and norm convergence of $\mcH^\theta(X, \overline{\mu})$}
It is well known that norm convergence implies modular convergence in classical Orlicz spaces (see \cite[Page 9]{Mali}). In this section, we discuss the relationship of modular and norm convergent of $\mcH^\theta(X, \overline{\mu})$ and $\mathcal{L}^\theta(X, \overline{\mu})$
\begin{Def}\label{def1}
We say that a sequence $(\xi_n) \in \mcH^\theta(X, \overline{\mu})$ is modular convergent to $\xi \in \mcH^\theta(X, \overline{\mu})$ if there exists a constant $k>0$ such that $\rho\big(k(\xi_n -\xi)\big) \to 0$ as $n \to \infty.$
\end{Def}
 
\begin{thm}\label{4th}
Modular convergent in $\mathcal{L}^\theta(X, \overline{\mu})$ implies modular convergent in $\mcH^\theta(X, \overline{\mu})$
\end{thm}
\begin{proof}
Let $\xi_n(.) \in \mathcal{L}^\theta( X,  \overline{\mu})$ be modular convergent to $\xi(.) \in \mathcal{L}^\theta(X, \overline{\mu}).$  Then there exists a constant $k >0$ such that $\lim\limits_{n \to \infty}\rho[ k(\xi_n - \xi)] =0.$ Since $\mathcal{L}^\theta(X, \overline{\mu}) $ is subset of $\mcH^\theta(X, \overline{\mu})$ as continuous dense embedding so, $\xi_n(.) \in H^\theta(X, \overline{\mu}).$ Using the definition \ref{def1},  $\xi_n(.)$ is modular convergent to $\xi(.) \in \mcH^\theta(X, \overline{\mu}).$
\end{proof}
\begin{rem}
The known fact norm convergence is modular convergence in $\mcL^\theta(X, \overline{\mu})$ along with the ( Theorem \ref{4th}), norm convergence in $\mcL^\theta(X, \overline{\mu})$ is modular convergence in $\mcH^\theta(X, \overline{\mu}).$
\end{rem}
Now we will check the relationship of norm convergence 
\begin{thm}
Norm convergent in $\mcL^\theta(X, \overline{\mu})$ are norm convergent in $\mcH^\theta(X, \overline{\mu}).$
\end{thm}
\begin{proof}
Let $\xi_n \in \mcL^\theta(X, \overline{\mu})$ such that $\xi_n \to f$ in $\mcL^\theta(X, \overline{\mu})$ in the way that $||\xi_n -\xi||_{\mcL^\theta} \to 0$ as $n \to \infty.$\\ This means,
\begin{align*}
&\inf\bigg\{k>0:~(L)\int_{\mathcal{J}}\theta\bigg(t, (\xi_n-\xi)(t)\bigg)d \overline{\mu} \leq 1\bigg\} \to 0\\ & \implies ~(L)\int_{\mathcal{J}}\theta\big(t, (\xi_n-\xi)(t)\big) \to 0\\& ~i.e.,~\theta(t, (\xi_n-\xi)(t)=0~\overline{\mu}-a.e.~as~n \to \infty\\&\xi_n-\xi=0~a.e.~as~n \to \infty 
\end{align*}
As,$~\xi_n \in \mcL^\theta(X, \overline{\mu}),~$ this implies $\xi_n \in \mcH^\theta(X, \overline{\mu})$ and $\xi_n-\xi=0 ~as~~\overline{\mu}-a.e..$ Now using the Theorem \ref{4rd}, $$(H) \int_{\mathcal{J}}(\xi_n-\xi) d \overline{\mu}=0.$$ Lastly, from  Definition \ref{def,}(d)
\begin{eqnarray*}
(H)\int_{\mathcal{J}}\theta(t, (\xi_n-\xi)(t))=0.
\end{eqnarray*}
So, \begin{eqnarray*}
\inf\bigg\{k>0:~(H)\int_{\mathcal{J}}\theta(t, (\xi_n-\xi)(t))d \overline{\mu} \leq 1\bigg\} \to 0.
\end{eqnarray*}
Hence, $||\xi_n-\xi||_{(X, \overline{\mu})} \to 0$ as $n \to \infty.$
\end{proof}
\begin{thm}
Modular convergence in $\mcL^\theta(X, \overline{\mu})$ is norm convergence in $\mcH^\theta(X, \overline{\mu}).$
\end{thm}
\begin{proof}
Let $\xi_n \in \mcL^\theta(X, \overline{\mu})$ be modular convergent to $\xi \in \mcL^\theta(X, \overline{\mu}).$ Then for $k>0$ (depending on the sequence of function $(\xi_n)$) such that $\rho\big(k(\xi_n-\xi)\big) \to 0$ as $n \to \infty.$ That is, 
\begin{align*}
&(L)\int_{\mathcal{J}}\theta(t, k(\xi_n-\xi))d\overline{\mu} \to 0~as~n \to \infty\\& i.e.,~\lim\limits_{n \to \infty}(L)\int_{\mathcal{J}}\theta(t, k(\xi_n-\xi))d\overline{\mu}=0\\& i.e.,~\theta(t, k(\xi_n-\xi))=0~\overline{\mu}~a.e.~when~n \to \infty.\\& {\textit{~For}}~k>0,~\\& \inf\bigg\{\frac{1}{k}>0:~(H)\int_{\mathcal{J}}\theta\bigg(t, \frac{\xi_n-\xi}{\frac{1}{k}}\bigg)d \overline{\mu}\leq 1 \bigg\}\to 0~\overline{\mu}-a.e.~as ~n \to \infty\\&So, ~||\xi_n-\xi||_{(X, \overline{\mu})} \to 0~as~n \to \infty.
\end{align*}
This complete the proof.
\end{proof}
\section*{Conclusion}
In this article we have discussed Henstock-Kurzweil integrable function space with a new norm equivalent to Alexiecz norm. $HK(X, \overline{\mu})$ is Banach space with the new norm (\ref{bh}). H-Orlicz space has been discussed with Banach valued Henstock-Kurzweil integrable function. We  have established the relationship of modular convergence and norm convergence of functions with values in Banach spaces in H-Orlicz space. We find modular convergence in $\mcL^\theta(X, \overline{\mu})$ is norm convergence in $\mcH^\theta(X, \overline{\mu}).$\\ We conclude this article with an open problem as follows:\\
{\bf{Problem:}}
Modular convergent in  $\mcH^\theta(X, \overline{\mu})$ does not implies modular convergent in  $\mathcal{L}^\theta(X, \overline{\mu}).$

\section{Declaration}
\noindent  {\bf Funding:} Not Applicable, the research is not supported by any funding agency.\\
 {\bf Conflict of Interest/Competing interests:} The authors declare that there is no conflicts of interest.\\
 {\bf Availability of data and material:} The article does not contain any data for
 analysis.\\
 {\bf Code Availability:} Not Applicable.\\
 {\bf Author's Contributions:} All the authors have equal contribution for the preparation of the article.


\end{document}